\documentclass[11pt]{amsart}

\usepackage{mathptmx}
\usepackage{multicol}

\usepackage{amsmath,amssymb,amsfonts,enumerate,amsthm,amscd,latexsym}

\newcommand{\Ann}{\mbox{Ann}\,}

\newcommand{\diam}{\mbox{diam}\,}
\newcommand{\girth}{\mbox{girth}\,}

\renewcommand{\dim}{\mbox{dim}\,}

\newcommand{\Z}{\mbox{Z}}

\renewcommand{\d}{\mbox{d}}

\newtheorem{thm}{Theorem}[section]
\newtheorem{cor}[thm]{Corollary}
\newtheorem{lem}[thm]{Lemma}
\newtheorem{prop}[thm]{Proposition}

\newtheorem{exam}[thm]{Example}
\newtheorem{rem}[thm]{Remark}

\begin{document}

\bibliographystyle{amsplain}

\author{Hamid Reza Maimani}
\address{Hamid Reza Maimani\\Department of Mathematics, University of
Tehran, Tehran, Iran\\ and Institute for Theoretical Physics and
Mathematics (IPM).}

\email{maimani@ipm.ir}

\author{Siamak Yassemi*}
\address{Siamak Yassemi\\Department of Mathematics, University of
Tehran, Tehran, Iran\\ and Institute for Theoretical Physics and
Mathematics (IPM).}

\email{yassemi@ipm.ir}


\thanks{* Corresponding author. Department of Mathematics,
University of Tehran, P.O. Box 13145--448 Tehran, Iran}



\title{Zero-divisor graphs of amalgamated duplication of a ring along an ideal}

\begin{abstract}
Let $R$ be a commutative ring with identity and let $I$ be an ideal
of $R$. Let $R\Join I$ be the subring of $R\times R$ consisting of
the elements $(r,r+i)$ for $r\in R$ and $i\in I$. We study the
diameter and girth of the zero-divisor graph of the ring $R\Join I$.

\end{abstract}

\maketitle

\section{Introduction}

Let $R$ be a commutative ring with non-zero unity. The concept of
the graph of the zero divisors of $R$ was first introduced by Beck
\cite{B}, where he was mainly interested in colorings. In his work
all elements of the ring were vertices of the graph. This
investigation of colorings of a commutative ring was then continued
by D.~D.~Anderson and Naseer in \cite{AN}. Let $\Z(R)$ be the set of
zero-divisors of $R$. In \cite{AL}, D.~F.~Anderson and Livingston
associate a graph, $\Gamma(R)$, to $R$ with vertices
$\Z(R)\setminus\{0\}$, the set of non-zero zero-divisors of $R$, and
for distinct $x,y\in\Z(R)\setminus\{0\}$, the vertices $x$ and $y$
are adjacent if and only if $xy=0$. Recall that a graph is said to
be {\it connected} if for each pair of distinct vertices $v$ and
$w$, there is a finite sequence of distinct vertices
$v=v_1,\cdots,v_n=w$ such that each pair $\{v_i,v_{i+1}\}$ is an
edge. Such a sequence is said to be a path and the distance,
$\d(v,w)$, between connected vertices $v$ and $w$ is the length of
the shortest path connecting them. The {\it diameter} of a connected
graph is the supremum of the distances between vertices. The
diameter is 0 if the graph consists of a single vertex, and a
connected graph with more than one vertex has diameter 1 if and only
if it is complete; i.e., each pair of distinct vertices forms an
edge. In \cite{AL}, the authors proved that $\Gamma(R)$ is always
connected and its diameter, $\diam(\Gamma(R))$, is always less than
or equal to 3 \cite[Theorem 2.3]{AL}. They also proved that
$\Gamma(R)$ is a complete graph if and only if either $R$ is
isomorphic to $\mathbb Z_2\times \mathbb Z_2$ or $xy = 0$ for all
$x,y\in\Z(R)$, cf. \cite[Theorem 2.8]{AL}. More recently, Axtell,
Coykendall and Stickles \cite{ACS} and Lucas \cite{L} have studied
the diameter of the corresponding graphs of the polynomial ring
$R[x]$ and the power series ring $R[[x]]$. Recall that the {\it
girth} of $G$ is the length of a shortest cycle in $G$ and is
denoted by $\girth(G)$. If $G$ has no cycles, we define the girth of
$G$ to be infinite. In \cite{AL}, the authors proved that the girth
of $\Gamma(R)$ is either infinite or less than or equal to four when
$R$ is Artinian and conjectured that this would hold if $R$ was not
Artinian, cf. \cite[Theorem 2.4]{AL}. In \cite[Theorem 1.6]{DS}
DeMeyer and Schneider and in \cite[Theorem 1.4]{M} Mulay proved this
conjecture independently (see also \cite[Theorem 2]{ACS}).

Let $M$ be an $R$-module, the idealization $R(+)M$ (also called
trivial extension), introduced by Nagata in 1956, cf, \cite{N}, is a
ring where the module $M$ can be viewed as an ideal such that its
square is $(0)$. In \cite{AS}, Axtell and Stickles considered zero
divisor graphs of idealization of commutative rings. They
characterize the diameter and the girth of the zero-divisor graph of
an idealization and show when this graph is complete.

In this paper, we deal with some applications of a similar general
construction, introduced recently in \cite{DF}, called the {\it
amalgamated duplication} of a ring $R$ along an ideal $I$, and
denoted by $R\Join I$. When $I^2=0$, the new construction $R\Join I$
coincides with the Nagata's idealization $R(+)I$. More precisely,
the amalgamated duplication of $R$ along an ideal $I$ is a ring that
is defined as the following subring of $R\times R$,
$$R\Join I=\{(r,r+i)|r\in R, i\in I\}.$$ More generally, this
construction can be given starting with a ring $R$ and an ideal $I$
of an overring $S$ of $R$ (such that $S\subseteq Q(R)$, where $Q(R)$
is the total ring of fractions of $R$); this extension has been
studied, in the general case, and form the different point of view
of pullbacks, by D'Anna and Fontana \cite{DF}. One main difference
of this construction, with respect to the idealization is that the
ring $R\Join I$ can be a reduced ring (and it is always reduced if
$R$ is a domain). As it happens for the idealization, one
interesting application of this construction is the fact that it
allows to produce rings satisfying (or not satisfying) preassigned
conditions. Moreover, in many cases, the amalgamated duplication of
a ring preserves the property of being reduced (see \cite{D},
\cite{DF}). Recently, D'Anna proved that, if $R$ is a Cohen-Macaulay
local ring, then $R\Join I$ is Gorenstein if and only if $I$ is a
canonical ideal, cf. \cite{D}, where was known for trivial
extension, cf. \cite{R}. This was our motivation to study the
zero-divisor graph of $R\Join I$.

In this paper we study the diameter and girth of the graph of
$\Gamma(R\Join I)$. In section 2, we review some properties of the
ring $R\Join I$ and classify the zero-divisors of this ring. In
section 3, we completely characterize the girth of zero-divisor
graph $\Gamma(R\Join I)$. More precisely, it is shown that $R$ is
not integral domain if and only if $\girth(\Gamma(R\Join I))=3$.
Also if $R$ is integral domain then $\girth(\Gamma(R\Join I))=4$
provided $|I|\ge 3$, and $\girth(\Gamma(R\Join I))$ is infinite if
$|I|=2$. In section 4, it is shown that for any non-zero ideal $I$
the following are equivalent:

\begin{itemize}

\item[(a)] The graph $\Gamma(R\Join I)$ is a complete graph;

\item[(b)] $(\Z(R))^2=0$ and $I\subseteq\Z(R)$;

\item[(c)] $(Z(R\Join I))^2=0$.

\end{itemize}

\section{Zero-divisors of the ring $R\Join I$}

Let $R$ be a commutative ring with identity element $1$ and let $I$
be an ideal of $R$. We define $R\Join I=\{(r,s)|r,s\in R, s-r\in
I\}$. It is easy to check that $R\Join I$ is a subring, with unit
element $(1,1)$, of $R\times R$ (with the usual componentwise
operations) and that $R\Join I=\{(r,r+i)|r\in R, i\in I\}$.

We recall that the idealization $R(+)M$, introduced by Nagata
\cite{N} for every $R$-module $M$, is defined as the $R$-module
$R\oplus M$ with multiplication defined by $(r,m)(s,n)=(rs,rn+sm)$.

It is easy to see that, if $\pi_i$ ($i=1,2)$) are the projections of
$R\times R$ on $R$, then $\pi_i(R\Join I)=R$ and hence if
$O_i=\ker(\pi_i|_{R\Join I})$, then $(R\Join I)/O_i\cong R$.
Moreover $O_1=\{(0,i)|i\in I\}$, $O_2=\{(i,0)|i\in I\}$ and $O_1\cap
O_2=(0)$. Now we state some properties of the ring $R\Join I$ from
\cite{DF}, that will be considered numerous times.

\begin{prop} (see \cite{DF}) The following hold:

\begin{itemize}

\item[(a)] The ring $R\Join I$ is reduced if and only if $R$ is reduced. In particular, if $R$ is an integral domain, $R\Join I$ is reduced and it has exactly two minimal primes which are $O_1$ and $O_2$.

\item[(b)] The ring $R\Join I$ is isomorphic to the idealization $R(+)I$ if and only if $I$ is a nilpotent ideal of index 2 in $R$.

\item[(c)] If in the $R$-module direct sum $R\oplus I$ we consider a multiplicative structure by setting $$(r,i)(s,j)=(rs,rj+si+ij),$$ then the map $f:R\oplus I\to R\Join I$ defined by $f((r,i))=(r,r+i)$ is a ring isomorphism.
So if we consider the ring $R\Join I$ as $R\oplus I$, and
$(r,i)(s,j)=(rs,rj+si+ij)$, then $O_1=\{(0,i)|i\in I\}$ and
$O_2=\{(-i,i)|i\in I\}$.
\end{itemize}

\end{prop}

In the rest of this paper we will use freely Proposition 2.1 part c
when we refer to the amalgamated duplication of $R$ along $I$.

To consider the zero-divisor graph of $R\Join I$ we need the
following result.

\begin{prop}

Let $R$ be a commutative ring and let $I$ be an ideal of $R$. Then
\[ \begin{array}{rl}
\Z(R\Join I)&=\{(0,i)|i\in I\}\cup\{(i,-i)|i\in I\}\\
&\cup\{(x,i)|x\in\Z(R)\setminus\{0\}, i\in I\}\\
&\cup\{(x,i)|x\in R\setminus\Z(R), \,\, \mbox{there exists} \,\, j\in I, j(x+i)=0\}.
\end{array} \]

\end{prop}

\begin{proof}

It is easy to see that $(0,i)$ is adjacent to $(j,-j)$ for any
$i,j\in I$. If $(a,i)\in R\Join I$ with $a\in\Z(R)$, then there
exists $b\neq 0$ such that $ab=0$. Consider the following two cases,
namely the case where $b$ does not belong to $\Ann(I)$ and the case
where $b$ belongs to $\Ann(I)$.

{\it Case 1.} Assume that $b\notin\Ann(I)$. There exists $j\in I$
such that $bj\neq 0$. We obtain $(a,i)(bj,-bj)=0$ and so
$(a,i)\in\Z(R\Join I)$.

{\it Case 2.} Assume $b\in\Ann(I)$. Then $bi=0$ and so
$(a,i)(b,0)=0$. Therefore $(a,i)\in\Z(R\Join I)$.

On the other hand, if $(a,i)\in\Z(R\Join I)$, $a\neq 0$ and
$a\in\Z(R)$, then there exists a non-zero element $s\in R$ such that
$as=0$. Thus $(a,i)\in\{(x,i)|x\in\Z(R)\setminus\{0\}, i\in I\}$. If
$a\notin\Z(R)$, then $(x,i)(y,j)=0$, implies that $y=0$ and
$j(i+x)=0$.

\end{proof}

\begin{rem}

Consider the following subsets of $\Gamma(R\Join I)$:

\begin{itemize}

\item[(1)] $T_1=\{(o,i)|i\in I\};$

\item[(2)] $T_2=\{(i,-i)|i\in I\};$

\item[(3)] $T_3=\{(x,i)|x\in\Z(R)\setminus\{0\}, i\in I\};$

\item[(4)] $T_4=\{(x,i)|x\in R\setminus\Z(R), j(x+i)=0\,\, \mbox{for some}\,\, j\in I\}.$

\end{itemize}

Then the following hold:

\begin{itemize}

\item[(a)] Each element of $T_1\setminus\{(0,0)\}$ is adjacent to any element of $T_2\setminus\{(0,0)\}$. This implies that there exists a complete bipartite graph $K_{|I|-1,|I|-1}$ in the structure of $\Gamma(R\Join I)$.

\item[(b)] If $i\in I\setminus\Z(R)$, then the vertices $(0,i)$ is adjacent only to vertices of $T_2\setminus\{(0,0)\}$, and $(i,-i)$ is adjacent only to vertices $T_1\setminus\{(0,0)\}$.

\item[(c)] There exists a subgraph of $\Gamma(R\Join I)$ isomorphic to $\Gamma(R)$.

\end{itemize}

\end{rem}

\section{Girth of $\Gamma(R\Join I)$}

In this section we study the girth of $\Gamma(R\Join I)$. If
$|I|=1$, then $\Gamma(R)=\Gamma(R\Join I)$ and so
$\girth(\Gamma(R))=\girth(\Gamma(R\Join I))$. Thus we are interested
in $\girth(\Gamma(R\Join I))$ for $|I|\ge 2$. The first result gives
complete answer for the rings that are not integral domain.

\begin{prop}

Let $I$ be an ideal of $R$. Then $\girth(\Gamma(R\Join I))=3$ if $R$
is not integral domain.

\end{prop}

\begin{proof}

Clearly $\Z(R)\cap I\neq\{0\}$. Consider $0\neq x\in\Z(R)\cap I$.
Then there exists $0\neq y\in R$ such that $xy=0$. Thus
$(0,x)$---$(y,0)$---$(x,-x)$---$(0,x)$ is a cycle of length 3 in the
graph $\Gamma(R\Join I)$. Therefore $\girth(\Gamma(R\Join I))=3$.

\end{proof}

\begin{prop}

Let $R$ be an integral domain and let $I$ be an ideal of $R$. Then
$\girth(\Gamma(R\Join I)=4$ provided $|I|\ge 3$. In addition, if
$|I|=2$, then $I=R\cong\mathbb Z_2$, and $\girth(\Gamma(R\Join
I))=\infty$.

\end{prop}

\begin{proof}

By assumption the only vertices of $\Gamma(R\Join I)$ are
$$\{(0,i)|i\in I\setminus\{0\}\}\cup\{(i,-i)|i\in
I\setminus\{0\}\}.$$ Thus $\Gamma(R\Join I)$ is a complete bipartite
graph. If $|I|\ge 3$, then for two distinct non-zero elements
$i,j\in I$, we have a cycle
$(0,i)$---$(i,-i)$---$(0,j)$---$(j,-j)$---$(0,i)$, in $\Gamma(R\Join
I)$ and hence $\girth(R\Join I)=4$.

Assume $|I|=2$. Then the graph $\Gamma(R\Join I)$ is isomorphic to
$(0,i)$---$(i,-i)$, and so $\girth(\Gamma(R\Join I)=\infty$.

\end{proof}

We obtain the following result by considering Propositions 3.1 and
3.2.

\begin{cor}
Let $I$ be a non-zero ideal of $R$. Then the following hold:

\begin{itemize}

\item[(a)] $\girth(\Gamma(R\Join I))=3$ if and only if $R$ is not an integral domain.

\item[(b)] $\girth(\Gamma(R\Join I))=4$ if and only if $R$ is an integral domain and $|I|\ge 3$.

\item[(c)] $\girth(\Gamma(R\Join I))=\infty$ if and only if $I=R\cong \mathbb Z_2$.

\end{itemize}

\end{cor}

\begin{cor}

The following statements are equivalent:

\begin{itemize}

\item[(a)] $R$ is integral domain.

\item[(b)] $\girth(\Gamma(R))=4$ or $\infty$.

\item[(c)] $R\Join I$ has exactly two minimal prime ideals $Q_1$ and $Q_2$ such that $Q_1\cap Q_2=(0)$ (i.e $Q_1=\{(0,i)|i\in I\}$ and $Q_2=\{(i,-i)|i\in I\}$).

\item[(d)] $\Gamma(R\Join I)$ is a complete bipartite graph.

\end{itemize}

\end{cor}

\begin{proof}

(a)$\Leftrightarrow$(b). This follows from Corollary 3.3.

(b)$\Rightarrow$(c). This follows from \cite[Proposition 2]{D}.

(c)$\Rightarrow$(d). This follows from \cite[Theorem 2.4]{AMY}.

(d)$\Rightarrow$(a). Since $\Gamma(R\Join I)$ is complete bipartite
graph, we have that $\girth(\Gamma(R\Join I)=4$ or $\infty$, so $R$
is integral domain.

\end{proof}

\section{Diameter of $\Gamma(R\Join I)$}

In this section we study the diameter of $\Gamma(R\Join I)$. It is
clear that if $\diam(\Gamma(R))>1$, then $\diam(\Gamma(R\Join
I))>1$. However, it is possible to have a ring such that
$\diam(\Gamma(R))=1$ but $\diam(\Gamma(R\Join I))\neq 1$, as the
following examples show.

\begin{exam}

Let $R=\mathbb Z_2\times\mathbb Z_2$ and $I=\mathbb Z_2\times\{0\}$.
Set $x=(1,0)$, $y=(0,1)$ and $z=(1,1)$. In the graph $\Gamma(R\Join
I)$, $(z,x)$ is adjacent exactly to $(0,x)$. On the other hand,
$(0,x)$ and $(x,0)$ are not adjacent. Thus $\diam(\Gamma(R\Join
I))=3$ but it is easy to see that $\diam(\Gamma(R))=1$.

\end{exam}

\begin{exam}
Let $R=\mathbb Z_6$ and $I=\{0,3\}$. Then $\diam(\Gamma(R))=2$. On
the other hand, $(1,3)\in\Z(R\Join I)$ is adjacent exactly to the
vertex $(0,3)$. Since $\d((0,3),(3,0))=2$, so $\d((1,3),(3,0))=3$
and hence $\diam(\Gamma(R\Join I))=3$.

\end{exam}

\begin{exam}
Let $R=\mathbb Z_8$ and $I=\{0,4\}$. Then $$\Z(R\Join
I)=\{(0,4),(4,4),(6,0),(2,0),(4,0),(2,4),(6,4)\}.$$ It is clear that
$\diam(R)=2=\diam(R\Join I)$.

\end{exam}

\begin{exam}

Let $R=\mathbb Z_2\times F$, where $F$ is a field. Let
$I=\{(0,0),(1,0)\}$. Then  $\Gamma(R)$ is a star graph by
\cite[Theorem 2.13]{AL}, and so $\diam(\Gamma (R))=2$. Consider the
element $((0,1),(1,0))\in\Z(R\Join I)$. It is clear that
$((0,1),(1,0))$ is adjacent to $((1,0),(1,0))$, and $((0,1),(1,0))$
is not adjacent to $((0,0),(1,0))$. Thus $((0,1),(1,0))$ is just
adjacent to $((1,0),(1,0))$. Now since $((1,0),(1,0))$ is not
adjacent to $((1,0),(0,0))$, so the distance of $((0,1),(1,0))$ to
$((1,0),(0,0))$ is equal 3. Therefore $\diam(\Gamma(R\Join I))=3$.

\end{exam}

\begin{lem}
Let $R$ be a commutative ring. Then $(Z(R\Join I))^2=0$ if and only
if $(\Z(R))^2=0$ and $I\subseteq\Z(R)$.

\end{lem}

\begin{proof}

``only if'' Let $x,y\in\Z(R)$. Then $(x,0), (y,0)\in\Z(R\Join I)$.
Thus $(x,0)(y,0)=(0,0)$ which implies $xy=0$. Therefore
$(\Z(R))^2=0$. Now assume $i\in I$. Then for any $z\in\Z(R)$, we
have $(0,i)(0,i)=(0,0)$. Thus $i^2=0$ and so $i\in\Z(R)$.

``if'' Suppose that $(\Z(R))^2=0$, and $I\subseteq\Z(R)$. For any
$x\in\Z(R)$, $i\in I$, the elements $(x,i)$ are adjacent to one
another. The only elements that we should study, are $(x,i)$ where
$x\in R\setminus\Z(R)$. Let $(x,i)\in\Z(R\Join I)$ and $x\in
R\setminus\Z(R)$. Then there exists $0\neq k\in I$ such that
$k(x+i)=0$. Since $I\subseteq\Z(R)$ and $(\Z(R))^2=0$, we have that
$kx=0$. Thus $x\in\Z(R)$, which is a contradiction. Therefore the
assertion holds.

\end{proof}

\begin{rem}
Note that in example 4.1, $\Gamma(R)$ is complete graph and
$I\subseteq\Z(R)$, but $\Gamma(R\Join I)$ is not complete.

\end{rem}

In the following Example, it is shown that the condition
$I\subset\Z(R)$ in Lemma 4.5 can not omitted.

\begin{exam}
Let $R=\mathbb Z_{p^2}$ where $p$ be a prime integer. It is easy to
see that $\Z(R)=\{0,p,2p,\cdots,p(p-1)\}$ and $(\Z(R))^2=0$. Let
$I=\mathbb Z_{p^2}$. then $(1,p-1)\in\Z(R\Join I)$ which is not
adjacent to $(p,0)$. Thus $(\Z(R\Join I))^2\neq0$.

\end{exam}

\begin{thm}

Let $I$ be a non-zero ideal of $R$. Then the following are equivalent:

\begin{itemize}

\item[(a)] The graph $\Gamma(R\Join I)$ is a complete graph

\item[(b)] $(\Z(R))^2=0$ and $I\subseteq\Z(R)$.

\item[(c)] $(Z(R\Join I))^2=0$.

\end{itemize}

\end{thm}

\begin{proof} (a)$\Rightarrow$(b). Assume $\Gamma(R\Join I)$ is a complete graph.
 Then $\Gamma(R)$ is a complete graph and so by \cite[Theorem 2.8]{AL} $R\cong\mathbb Z_2\times\mathbb Z_2$ or $xy=0$
  for all $x,y\in\Z(R)$. If $R\cong\mathbb Z_2\times\mathbb Z_2$ then for any non-zero ideal $I$ of $\mathbb Z_2\times\mathbb Z_2$, two vertices $((0,1),(0,0))$ and $((0,0),(0,1))$ are not adjacent
 and hence $\Gamma(R\Join I)$ is not complete graph. Thus $xy=0$ for all $x,y\in\Z(R)$, that means $(\Z(R))^2=0$.
 Observing that for any $i,j \in I, (0,i)(0,j)=0$ and hence $ij=0$
 which implies $I\subseteq Z(R)$

(b)$\Rightarrow$(c).  The assertion follows from Lemma 4.5.

(c)$\Rightarrow$(a). This is clear.

\end{proof}

\begin{lem}

Assume $R$ is not integral domain and $I\nsubseteq\Z(R)$. If $\Z(R)$
is an ideal then $\diam(\Gamma(R\Join I))=3$.

\end{lem}

\begin{proof}
Choose $i\in I\setminus\Z(R)$. Since $I\cap\Z(R)\neq (0)$, there exists $k\in I$, such that
 $Ann_R(k)\neq (0)$. Let $0\neq x\in\Ann_R(k)$. Then $k(x-i+i)=0$. Set $y=x-i$. The vertex $(y,i)$ is adjacent
to $(0,k)$. If $y\notin\Z(R)$, then $(y,i)$ is adjacent exactly to
vertices $(0,l)$ where $y+i\in\Ann_R(l)$.
 Since the distance of such vertices $(0,l)$ to $(0,i)$ is 2, we have $\diam(\Gamma(R\Join I))=3$.
 If $x-i=y\in\Z(R)$, then we have $i\in\Z(R)$, since $x\in\Z(R)$ - a contradiction.

\end{proof}

\begin{cor}
If $I\nsubseteq\Z(R)$ and there is a vertex of $\Gamma(R)$ which is adjacent to every other vertex of $\Gamma(R)$,
 then $\diam(\Gamma(R\Join I))=3$.

\end{cor}

\begin{proof}
Since there exists a vertex of $\Gamma(R)$ which is adjacent to
every other vertex, by \cite[Theorem 2.5]{AL}, $\Z(R)$ is an
annihilator ideal or $R\cong\mathbb Z_2\times A$, where $A$ is an
integral domain. If $\Z(R)$ is an ideal, then $\diam(\Gamma(R\Join
I))=3$ by Lemma 4.9. If $R\cong\mathbb Z_2\times A$, then
$I=\{0\}\times J$ or $I=Z_2\times J$ where $J$ is a non-zero ideal
of $A$. Since $\{0\}\times J\subseteq Z(R)$ and $J\nsubseteq Z(R)$,
we have $I=Z_2\times J$. If $A=\Bbb Z_2$, then $I=R\cong\mathbb
Z_2\times\mathbb Z_2$, and hence $\diam(\Gamma(R\Join I))=3$.

Now assume that $A\neq\Bbb Z_2$. Let $A\neq J$ and $a\in A\setminus
J$. For arbitrary non-zero element $b\in J$, consider $X=(1,a)$,
$Y=(1,b)$, $W=(b,b)$, and $Z=(1,0)$. Then $(X,Y)$ is adjacent
exactly to the vertex $(0,Z)$, and the distance between $(0,Z)$ and
$(0,W)$ is equal 2. Therefore, $\diam(\Gamma(R\Join I))=3$.

In the case $A=J$, consider $X=(1,a)$, $Y=(1,b)$, $W=(1,1)$, and
$Z=(1,0)$, where $a,b\in A\setminus\{0\}$. Then $\d((X,Y), (0,W))=3$
and hence $\diam(\Gamma(R\Join I))=3$.

\end{proof}

In \cite[Example 3.7]{AS}, the authors give an example of a ring $R$
and an $R$-module $K$ with $\diam(\Gamma(R))=3$ but
$\diam(\Gamma(R(+)K))=2$. In the following result we show that this
case does not happen for $\Gamma(R\Join I))$.

\begin{prop}

If $\diam(\Gamma(R))=3$, then $\dim(\Gamma(R\Join I))=3$.

\end{prop}

\begin{proof}

Let $\diam(\Gamma(R\Join I)<3$. Then clearly $\diam(\Gamma(R\Join
I))=2$.
 Choose $x,y\in\Z(R)$ with $d(x,y)=3$ (in the graph $\Gamma(R)$). Consider two vertices
 $(x,0), (y,0)$ in $\Gamma(R\Join I)$. Then $(x,0)$ and $(y,0)$ are not adjacent.
 Thus $d((x,0),(y,0))=2$, so there exists $(r,i)\in\Z(R\Join I)$ such that $(r,i)$ adjacent with two
  vertices $(x,0)$ and $(y,0)$. Hence $r,i\in\Ann_R(x)\cap\Ann_R(y)$. Since $r,i\notin\{x,y,0\}$,
  we have $d(x,y)\leq 2$ in $\Gamma(R)$, which is a contradiction.

\end{proof}

\begin{thm}
If $\Z(R)$ is not an ideal, then $\diam(\Gamma(R\Join I))=3$.

\end{thm}

\begin{proof} First suppose that $R$ is reduced ring. We follow the following two steps.

{\it Step 1.} We show that $\diam(\Gamma(R\Join I))\neq 2$. Let
$\diam(\Gamma(R\Join I))=2$. Since $R$ is reduced,
 $R\Join I$ is reduced, cf. \cite[Proposition 2]{D}. On the other hand, since $\Z(R)$ is not an ideal,
 there exists $x,y\in\Z(R)$ such that $x-y\notin\Z(R)$, and so $(x,0)-(y,0)\notin\Z(R\Join I)$.
 Hence $\Z(R\Join I)$ is not an ideal of the ring $R\Join I$. Now by \cite[Theorem 2.2]{L}, $R\Join I$
 has exactly two two minimal prime ideals $P_1,P_2$. Therefore, $P_1\cap P_2=\{0\}$, and $\Gamma(R\Join I)$
 is complete bipartite graph. Thus by Corollary 3.4, $R$ is an integral domain and so $\Z(R)$ is an ideal, which is a contradiction.

{\it Step 2.} We show that $\diam(\Gamma(R\Join I))\neq 1$. Let
$\diam(\Gamma(R\Join I))=1$. Then $R\Join I$ has exactly two minimal
ideals by \cite[Theorem 2.2]{L}, and so $R$ is integral domain by
Corollary 3.4. This is a contradiction.

Now suppose that $R$ is not reduced ring. Then $\diam(\Gamma(R))=3$
by \cite[Corollary 2.5]{L} and so $\diam(\Gamma(R\Join I))=3$ by
Proposition 4.11.

\end{proof}

\begin{prop}
Let $\Z(R)$ be an ideal of $R$ and $I\subseteq\Z(R)$. For all
adjacent vertices $a,b$ of $\Gamma(R)$, let $\Ann(a,b)\neq 0$. Then
$\diam(\Gamma(R\Join I))=2$ provided $\diam(\Gamma(R))=2$.

\end{prop}

\begin{proof}
Since $\diam(\Gamma(R))=2$, so $\diam(\Gamma(R\Join I))\ge 2$. Let
$(x,i)$ and $(y,j)$ be two vertices of $\Gamma(R\Join I)$. Consider
the following cases:

{\it Case 1.} Let $x=y=0$. Then $(0,i)$ and $(0,j)$ are adjacent to all vertices $(k,-k)$ where $k\in I$. Therefore $\d((0,i),(0,j))\le 2$.

{\it Case 2.} Let $x=0$ and $y\neq 0$. Since $(y,j)\in\Z(R\Join I)$,
we claim that $y\in\Z(R)$. If not, $(y,j)$ is adjacent to vertices
$(0,k)$ where $k(y+j)=0$. Thus $y+j\in\Z(R)$. Since
$I\subseteq\Z(R)$, we have that $y\in\Z(R)$ which is a
contradiction. Therefore $y\in\Z(R)$. If there exists a non-zero
element $z\in\Ann(y)$, and a non-zero element $k\in I$ such that
$zk\neq 0$, then we have path $(0,i)$---$(zk,-zk)$---$(y,j)$. If for
any $z\in\Ann(y)$ and $k\in I$ we have $zk=0$, then for an element
$0\neq z\in\Ann(y)$ we have path $(y,j)$---$(z,0)$---$(0,i)$.
Therefore $\d((y,j),(0,i))\le 2$.

{\it Case 3.} Let $x\neq 0$ and $y\neq 0$. By the same argument as
Case 2, $x,y\in\Z(R)$. If $\d(x,y)= 2$, then there exists $0\neq
z\in\Z(R)$ such that $yz=xz=0$. If there exists $k\in I$ such that
$zk\neq 0$, then $(x,i)$ and $(y,j)$ are adjacent to $(zk,-zk)$, and
hence $\d((x,i),(y,j))\le 2$. If $zk=0$ for each $k\in I$, then we
have path $(x,i)$---$(z,0)$---$(y,j)$, and the assertion holds.

If $\d(x,y)=1$, then $x$ and $y$ are adjacent in $\Gamma(R)$. Thus
there exists $0\neq z\in\Z(R)$ such that $xz=yz=0$, since
$\Ann(x,y)\neq 0$. So by a same argument as above,
$\d((x,i),(y,j))\le 2$. Therefore $\d(a,b)\le 2$ for any
$a,b\in\Gamma(R\Join I)$ and hence $\diam(\Gamma(R\Join I)=2$.

\end{proof}

\begin{cor}
Let $R$ be non-reduced ring, $\Z(R)$ be an ideal of $R$, and
$I\subseteq\Z(R)$. Then $\diam(\Gamma(R\Join I)=2$ provided
$\diam(\Gamma(R))= 2$.

\end{cor}

\begin{proof}
Let $a,b\in\Z(R)$. If $a$ and $b$ are adjacent and $\Ann(a,b)=0$, then by \cite[Theorem 2.4]{L} $\diam(\Gamma(R))=3$. This is a contradiction. Therefore for any two adjacent vertices $a$ and $b$ of $\Gamma(R)$ we have that $\Ann(a,b)\neq 0$. Now use Proposition 4.13.

\end{proof}

\begin{lem}

If $I\nsubseteq\Z(R)$ and $\diam(\Gamma(R\Join I))=2$, then for any
$y\in\Z(R)\setminus\{0\}$, $\Ann_R(y)\cap I\neq\{0\}$.

\end{lem}

\begin{proof}

Let $i\in I\setminus\Z(R)$ and $y\in\Z(R)\setminus\{0\}$. Then the
vertices $(0,i)$ and $(y,0)$ are not adjacent. Thus there exists
$(s,j)\in\Z(R\Join I)\setminus\{(0,0)\}$ such that $(s,j)$ adjacent
to both vertices. So $i(j+s)=0$. Since $i\notin\Z(R)$, we have that
$0\neq s=-j\in I$. In addition, $sy=0$ implies that $s\in\Ann_R(y)$.
Therefore $0\neq s\in\Ann_R(y)\cap I$.

\end{proof}

Our last result provides a condition which is sufficient for $\Z(R)$ to be a prime ideal (that means $R$ has exactly one associated prime).

\begin{prop}

If there exists an element $(r,i)\in R\Join I$ which is adjacent to
every vertices of $\Gamma(R\Join I)$, then $\Z(R)$ is a prime ideal.

\end{prop}

\begin{proof}

For any $x\in\Z(R)$, we have $(r,i)(x,0)=0$. Thus $ix=0$ for any $x\in\Z(R)$. If $i\neq 0$ then $\Z(R)=\Ann_R(i)$. If $i=0$, then $r\neq 0$ and $rx=0$ for any $x\in\Z(R)$. This means that $\Z(R)=\Ann_R(r)$. therefore $\Z(R)$ is an ideal. On the other hand $R\setminus\Z(R)$ is a multiplicative closed subset of $R$ and so $\Z(R)$ is prime ideal.

\end{proof}



\providecommand{\bysame}{\leavevmode\hbox
to3em{\hrulefill}\thinspace}


\begin{thebibliography}{10}

\bibitem{AL}
D.~F.~Anderson and P.~S.~Livingston, \emph{The Zero-divisor graph of
a commutative ring}, J. Algebra  \textbf{217}  (1999), 434--447.

\bibitem{AMY}
S.~Akbari, H.~R.~Maimani and S.~Yassemi, \emph{When a zero-divisor
graph is planar or a complete $r$-partite graph}, J. Algebra
\textbf{270} (2003), 169--180.

\bibitem{AN}
D.~D.~Anderson and M. Naseer, \emph{Beck's coloring of a commutative
ring}, J. Algebra \textbf{159} (1993), 500--514.

\bibitem{B}
I.~Beck, \emph{Coloring of Commutative Rings}, J. Algebra
\textbf{116} (1988), no. 1, 208-226.

\bibitem{ACS}
M.~Axtell, J.~Coykendall and J.~Stickles, \emph{Zero-divisor graphs
of polynomial and power series over commutative rings}, Comm.
Algebra \textbf{33} (2005), 2043--2050.

\bibitem{AS}
M.~Axtell and J.~Stickles, \emph{Zero-divisor graphs of
idealizations}, J. Pure Appl. Algebra, \textbf{204} (2006),
235--243.

\bibitem{D}
M.~D'Anna, \emph{A construction of Gorenstein rings}, J. Algebra
\textbf{306} (2006), 507--519.

\bibitem{DF}
M.~D'Anna and M.~Fontana, \emph{An amalgamated duplication of a ring
along an ideal}, to appear in J. Algebra Appl.

\bibitem{DS}
F.~DeMeyer and K. Schneider, \emph{Automorphisms and zero-divisor
graphs of commutative rings}, Int. J. Commut. Rings \textbf{1}
(2002), 93--106.

\bibitem{L}
T.~G.~Lucas, \emph{The diameter of a zero divisor graph}, J. Algebra
\textbf{301} (2006), 174--193.

\bibitem{M}
S.~B.~Mulay, \emph{Cycles and symmetries of zero-divisors}, Comm. Algebra \textbf{30} (2002), 3533--3558.

\bibitem{N}
M.~Nagata, \emph{Local rings, Interscience}, New York, 1962.

\bibitem{R}
 I.~Reiten, \emph{The converse to a theorem of Sharp on Gorenstein modules},  Proc. Amer. Math. Soc. \textbf{32} (1972), 417--420.

\end{thebibliography}
\end{document}